\documentclass{amsart}
\usepackage{amssymb, amsmath}

\newtheorem{thm}{Theorem}[section]
\newtheorem{lemma}[thm]{Lemma}

\newtheorem{prop}[thm]{Proposition}
\newtheorem{cor}[thm]{Corollary}
\newtheorem*{assumpt}{Assumption}

\theoremstyle{definition}

\renewcommand{\Re}{\operatorname{\rm Re}\nolimits}
\renewcommand{\Im}{\operatorname{\rm Im}\nolimits}

\def \comp {\operatorname{comp}}
\def \supp {\operatorname{supp}}
\def \diam {\operatorname{diam}}
\def \tr {\operatorname{tr}}

\def \restrict {\upharpoonright}
\def \Real {{\mathbb R}}
\def \Sphere {\mathbb{S}}
\def \Complex {\mathbb{C}}
\def \Natural {{\mathbb N}}
\def \Sphere {{\mathbb S}}
\def \complement {{\bf C}}
\def \mcr{{\mathcal R}}
\def \hol {{\mathcal H}}
\def \Integers {{\mathbb Z}}\title [Resonances in potential scattering]
{Several complex variables and the distribution of resonances in 
potential scattering}
   \author { T. Christiansen}
\thanks{Partially supported by NSF grant DMS 0088922.}
\thanks{e-mail address: {\tt tjc@math.missouri.edu}}
\begin{document}
\begin{abstract} 
We study resonances associated to Schr\"odinger operators with compactly
supported potentials on $\Real^d$, $d\geq3$, odd.  
We consider compactly supported potentials depending holomorphically on
a parameter
$z\in \Complex^m$.  For certain such families, for all 
$z$ except those in a pluripolar set, the associated 
resonance-counting function has order of growth $d$.
\end{abstract}

\maketitle

\section{Introduction}

In this paper we study the growth of the resonance-counting function
for potential scattering in odd dimension $d\geq 3$.
Let $V\in L^{\infty}_{\comp}(\Real^d;\Complex)$ and let $N_V(r)$ be
the resonance-counting function for the Schr\"odinger operator $\Delta +V$.
The purpose of this paper is to show that
\begin{equation}\label{eq:biggrowth}
\lim \sup _{r\rightarrow \infty}\frac{\log N_V(r)}{\log r}=d
\end{equation}
for many potentials $V$.  By \cite{zwodd}, this is the maximum value this
limit can obtain.  Previously, the only potentials known to
satisfy (\ref{eq:biggrowth}) in
dimension at least three were a class of radial potentials \cite{zwrp}.
For a certain class of compactly supported
potentials $W(z)$ depending holomorphically on a parameter $z$, we show that
(\ref{eq:biggrowth}) holds for $V=W(z)$, for all $z$ except those in a 
pluripolar set.  In a probabilistic sense this greatly expands the number
of potentials which are known to have resonance-counting function
with maximal order of growth.    
We use this to show that potentials with this 
property are dense in the $L^{\infty}$ norm in 
$L^{\infty}_{\comp}(\Real^d)$.  We remark that there are complex-valued
$V\in L^{\infty}_{\comp}(\Real^d;\Complex)$ such that the limit in
(\ref{eq:biggrowth}) is $0$ \cite{sownr}.

For odd $d$ and $V\in L^{\infty}_{\comp}(\Real^d;\Complex), $
$R_V(\lambda)=(\Delta+V-\lambda^2)^{-1}$ is defined for $\Im \lambda >0$.
 If $\chi \in C_c^{\infty}(\Real^d)$, $\chi \equiv 1$ on the support of 
$V$, then $\chi R_V\chi$ has a meromorphic continuation to $\Complex$.  The
poles of this continuation are resonances, or scattering poles.  They
are, in many ways, analogous to eigenvalues and correspond to decaying 
states.  For an introduction to resonances and for a survey of some 
results on their distribution, see \cite{vodevsurvey, zwsurvey, zwnotices,
zwqr}.

Let $\mcr_V$ be the 
set of poles of $R_V(\lambda)$, repeated with multiplicity.
 Let
$$N_V(r)=\# \{ z_j\in \mcr_V:\; |z_j|<r\}.$$  Then, if $d=1$, 
$$\lim _{r\rightarrow \infty}\frac{N_V(r)}{r}=\frac{2}{\pi }\diam (\supp(V))
$$
\cite{froese, zw1,simon}.
This is true for complex-valued potentials as well as for real-valued ones.
Much less is known about the higher-dimensional case, and there is 
evidence that the question of distribution of resonances is more subtle.
Zworski \cite{zwodd} showed that for $d\geq 3$, odd, 
$$N_V(r)\leq C (r^d+1)$$
and this order of growth is achieved by a class of radial potentials 
\cite{zwrp}.  On the other hand, the best known lower bound to hold for
a general class of potentials is, for non-trivial $V\in C^{\infty}_{c}
(\Real^d;\Real)$,
$$\lim\sup _{r\rightarrow \infty}\frac{N_V(r)}{r}>0$$
\cite{sBlb}.  It is important that these are {\em real}-valued
potentials, as this does not hold for all smooth complex-valued potentials.  
In \cite{sownr}, there is an example of a family of {\em complex}-valued
potentials for which $N_V(r)\equiv 0$ for all $r$.  
In fact, the example works in
even dimensions as well (with some caveats for 
$d=2$).  The potentials can be chosen 
to be smooth.

	In this paper we show that there are many potentials with
resonance-counting function with 
the maximum order of growth.
This theorem
 can be viewed as providing a kind of probabilistic lower bound on the
resonance counting function, as it gives no information for any given
potential but says that ``most'' potentials in certain families have 
resonance counting function with maximal growth rate.  
The proof uses some results from several 
complex variables.
\begin{thm}\label{thm:intermediate}
Let $d\geq 3$ be odd and 
let $\Omega \subset \Complex^m$ be an open, connected set.  Let
$V(z,x)=\sum_{j=1}^{j_0} f_j(z)V_j(x)$ with $f_j$ holomorphic on $\Omega$ and
$V_j\in L^{\infty}_{\comp}(\Real^d;\Complex)$.  Suppose
$$\lim \sup _{r \rightarrow \infty}\frac{\log N_{V(z_0)}(r)}{\log r} 
=d$$
for some $z_0 \in \Omega$.  Then
$$\lim \sup _{r \rightarrow \infty}\frac{\log N_{V(z)}(r)}{\log r} 
=d$$ for $z\in \Omega \setminus E$, where $E$ is a pluripolar set.
\end{thm}
We recall the definition of a pluripolar set in Section
\ref{s:sca} and refer
the reader to \cite{klimek, l-g} for further details.  We remark that
pluripolar
sets are quite small-- in particular, they have Lebesgue measure zero, 
and there are further restrictions on them.

We remark that by \cite[Theorem 2]{zwrp} the condition on
$V(z_0)$ is satisfied if $V(z_0)$ is
the radial potenial 
$W(|x|)$, with $W\in C^2[0,a]$, $W(a)\not =0$.  One may
also use Theorem \ref{thm:fixedsign} to generate such potentials.

As an application of Theorem \ref{thm:intermediate} and some further
study of holomorphic functions whose zeros correspond to resonances,
we obtain the following theorem.
\begin{thm} \label{thm:fixedsign} Suppose
$d$ is odd and $V\in L^{\infty}_{\comp}(\Real^d;\Real)$ is bounded 
below by the characteristic function of a ball.  Then
$$\lim \sup _{r\rightarrow \infty}\frac{\log N_{zV}(r)}{\log r}=d$$
for $z\in \Complex \setminus E$, where $E\subset \Complex$ is a pluripolar
set.
\end{thm}  Earlier results for potentials of fixed sign are found in 
 \cite{l-pdm} and \cite{vasy}.  These papers studied the purely 
imaginary
scattering poles associated to potentials 
$V\in L^{\infty}_{\comp}(\Real^d;\Real)$ where $V$ or $-V$ is bounded below by a positive multiple of the 
characteristic function of a ball.  They showed that for such potentials,
$$\# \{ \lambda_j\in \mcr_{V}:\; \lambda_j\in i\Real,\; |\lambda_j|\leq r\}
\geq c_Vr^{d-1}$$
for some constant $c_V>0$.

A corollary of Theorem \ref{thm:intermediate}, Theorem \ref{thm:fixedsign},
and the properties of pluripolar sets is
\begin{cor} \label{t:density}
 For $d\geq 3$, odd, 
the set 
$$\{ V\in L^{\infty}_{\comp}(\Real^d;\Real):\; \lim\sup_{r\rightarrow \infty}
\frac{\log N_V(r)}{\log r} =d \}$$
is dense in $L^{\infty}_{\comp}(\Real^n;\Real)$ under
the $L^{\infty}$ norm.  The set
$$\{ V\in L^{\infty}_{\comp}(\Real^d;\Complex):\; \lim\sup_{r\rightarrow \infty}
\frac{\log N_V(r)}{\log r} =d \}$$
is dense in the set $ L^{\infty}_{\comp}(\Real^d;\Complex)$ under
the $L^{\infty}$ norm.  Moreover, the same results are true if we replace
$L^{\infty}_{\comp}$ by $C_c^{\infty}$ and the $L^{\infty}$ norm by the
$C^{\infty}$ topology.
\end{cor}
We remark that to prove the results for $L^{\infty}_{\comp}$ potentials in the
$L^{\infty}$ topology, one could use
\cite[Theorem 2]{zwrp} instead of Theorem \ref{thm:fixedsign}.

 In the next section of this paper,
 we recall some definitions and
facts from one and 
several complex variables.  In addition, we prove an extension
of \cite[Corollary 1.42]{l-g}, a 
result about order of growth for functions of several complex
variables.  This result, combined with some facts
about the determinant of the scattering matrix which are established in
Section \ref{s:sd}, enables us to prove our main theorem,
Theorem \ref{t:major1}, in Section 
\ref{s:pot}.  This result is somewhat stronger than Theorem 
\ref{thm:intermediate}.  Section \ref{s:fixedsign}
is devoted to the proofs of Theorem \ref{thm:fixedsign} and Corollary
\ref{t:density}.

Throughout this paper, $C$, $C_V$, $C_{\epsilon}$, and $C_{\chi}$ denote
constants whose value may change from line to line.  The dimension $d$
is odd throughout. 

We are pleased to thank D. Drasin, D. Edidin,
 C. Kiselman, and I. Verbitsky for 
helpful discussions.

\section{Some complex analysis}\label{s:sca}

In this section we recall some definitions and results from complex analysis,
and prove an extension of a result in several complex variables that we
shall need.  

Let $a_1, a_2, a_3....$ be a sequence of non-zero
complex numbers with $|a_m|\rightarrow
\infty$.  The {\em convergence exponent} of this sequence is the 
greatest lower bound of the set
$$\left\{ \lambda: \sum_{m=1}^{\infty}\frac{1}{|a_m|^\lambda}\; \text{
converges}\right\}.$$ If $n(r)=\#\{ a_j: |a_j|<r\}$,
then $$\limsup _{r\rightarrow \infty}\frac{\log n(r)}{\log r},$$
which may be called the order of $n(r)$, 
is the same as the convergence exponent for the sequence 
$\{ a_j\}_{j=1}^{\infty}$.  We shall abuse notation slightly and call this 
the convergence exponent for the {\em set} $\{a_l\}$, where we order 
$\{a_l\}$ so that $|a_1|\leq |a_2| \leq |a_3|\leq ...$ to form the sequence.

We now recall the definition and
some facts about plurisubharmonic functions.
For further details, see, for example, \cite{klimek,l-g}.

Let $\Omega \subset \Complex^m$ be a domain; that is, an open, connected
set.  A function $\varphi(z)$ which
takes its values in $[-\infty, \infty)$ is {\em plurisubharmonic} in $\Omega$
if 
\begin{itemize}
\item
 $\varphi(z)$ is upper semi-continuous and $\varphi \not \equiv -\infty$.
\item 
 For every $z\in \Omega$ and every $r$ such that $\{ z+uw: |u|\leq r,\; 
u \in \Complex\}\subset \Omega$, 
$$\varphi(z)\leq 
(2\pi)^{-1} \int_0^{2\pi}\varphi(z+re^{i\theta}w)d\theta.$$
\end{itemize}  We shall write $\varphi \in \text{PSH}(\Omega)$ if
$\varphi$ is plurisubharmonic on $\Omega.$
Being plurisubharmonic is a local property.  
Let $\Omega \subset \Complex^m$ be a domain.  If 
$\varphi$ is upper semi-continuous on $\Omega$, $\varphi \not \equiv 
-\infty$, 
and for every $z\in \Omega$ there is a $\rho(z)$ such that 
$$\varphi(z) \leq (2\pi)^{-1}\int _0^{2\pi}\varphi(z+we^{i\theta})d\theta
$$
for all $w\in \Complex^m$, $\|w\| <\rho(z)$, then we say that 
$\varphi$ is locally plurisubharmonic on $\Omega$.  But if $\varphi$ is
locally plurisubharmonic on $\Omega$, it is plurisubharmonic on $\Omega$
(e.g. \cite[Proposition I.19]{l-g}).

A set $E\subset \Complex^m$ is {\em pluripolar} if for each $a\in E$
there is a neighborhood $V$ of $a$ and $\varphi \in \text{PSH}(V)$ 
such that $E\cap V\subset \{z\in V:\varphi(z)=-\infty\}$.  This is 
equivalent to the definition given in \cite{l-g} via the Josefson Theorem
\cite[Theorem 4.7.4]{klimek}.

For a function $\varphi$ which is plurisubharmonic in 
$\theta_1 < \arg u < \theta_2$, we define the order $\rho$ of $\varphi$ 
 in $\theta_1<\arg u<\theta_2$ as 
$$\rho =\lim \sup_{
r \rightarrow \infty}
 \frac{\log \sup_{ \theta_1<\arg u <\theta_2, |u|=r} |\varphi(u)|}{\log r}.$$
An important example of a plurisubharmonic function is $\log |f|$, where
$f$ is holomorphic.  Thus we shall make the following (standard) definition
of order for a holomorphic function.
For $f$ holomorphic for $\theta_1 < \arg u < \theta_2$, the order $\rho$ of
$f$ in $\theta_1 <\arg u <\theta_2$ is 
$$\rho =\lim \sup_{
r \rightarrow \infty}
 \frac{\log \sup_{ \theta_1<\arg u <\theta_2, |u|=r}\log |f(u)|}{\log r}.$$
Since the two notions of order are so closely related, we use the 
same name and notation for each.

We shall be concerned with functions that satisfy the following set of 
assumptions.
\begin{assumpt}[{\bf A0}]
For some 
open $\Omega \subset \Complex^m$ and some $\epsilon>0$,  
$f(z,\lambda)$ is holomorphic on 
$\Omega \times \{ \lambda \in \Complex: \Im \lambda >-\epsilon\}$.  
Moreover, there are constants $C_f$
and $\alpha$ such that
\begin{equation}\label{eq:realbound}
\log |f(z,\lambda)|\leq C_f(1+|\lambda|^{\alpha})\; \text{for}\;
 \lambda \in \Real.
\end{equation}
\end{assumpt}
With the next two lemmas, we construct a plurisubharmonic function 
on $\Omega \times \Complex$ whose
order is related to the order of the function $f$ in a half-plane.
Some related results and techniques appear in Theorem I.28 and its proof 
in \cite{l-g}.

\begin{lemma}\label{l:M}  Assume $f$ satisfies assumption (A0).  For some 
$\beta >\alpha$, $\beta\geq 1$, let
$$M(z,r)=\max (\max_{\genfrac{}{}{0pt}{}{|\lambda|\leq r }{ \Im \lambda \geq 0}}
\log |f(z,\lambda)|, 
r^{\beta}).$$  Then there is an $r_0\in \Real$ such that 
 $M(z,r)$ is a plurisubharmonic function of $(z,u)\in \Omega \times 
\{ u\in \Complex: |u|>r_0\}$, where $|u|=r$.
\end{lemma}
\begin{proof}
Note that since $f$ is holomorphic, we actually have that $M$ is
continuous.  Clearly, $M\not \equiv -\infty$.

Now we shall show that $M$ is locally plurisubharmonic when $r$ is 
sufficiently large.  The key idea here is that we require $r$ to be so
large that 
$$C_f(1+r^{\alpha})< r^{\beta}$$ 
where $C_f$, $\alpha$ are as in (\ref{eq:realbound}).
For such values of 
$r$, $M(z,r)$ is not the value of $\log|f(z,\pm r)|$.

Choose $r$ sufficiently large as above.
Suppose $M(z,r)=\log|f(z,\lambda_0)|$, with $|\lambda_0|=r$ and 
$\Im \lambda_0>0$.  Then to
see that $M$ is locally plurisubharmonic at $(z,u)$ with $|u|=r$, 
consider $w\in \Complex^{m+1}$ with $\|w\|\leq \Im \lambda_0$
and write $w=(w',w_{m+1})$.  Note that $(w', \lambda_0w_{m+1}/u)$
has the same norm as $w$ and 
$$|\lambda_0+ e^{i\theta}\lambda_0w_{m+1}/u|=
|u+e^{i\theta} w_{m+1}|.$$  Then
\begin{align*}
M(z,|u|)& = M(z,|\lambda_0|)\\ & 
\leq (2\pi)^{-1}\int_0^{2\pi}
 \log |f(z+e^{i\theta} w', \lambda_0+e^{i\theta}
\frac{ \lambda_0}{u}w_{m+1})|d\theta  \\
&  
\leq (2\pi)^{-1}\int_0^{2\pi} M(z+e^{i\theta} w', |\lambda_0+e^{i\theta}
\frac{ \lambda_0}{u}w_{m+1}|)d\theta \\
& = (2\pi)^{-1}\int_0^{2\pi} M(z+e^{i\theta} w', |u+e^{i\theta}
w_{m+1}|)d\theta.
\end{align*}

Suppose, on the other hand, that $M(z,r)=r^{\beta}$. 
Let $w=(w',w_{m+1})\in \Complex^{m+1}$.  Then, if 
$|u|=r$, 
\begin{align*}
M(z,|u|)& = M(z,r)=r^{\beta}\\
& \leq (2\pi)^{-1}\int _0^{2\pi} |u+w_{m+1}e^{i\theta}|^{\beta}d\theta\\
& \leq 
(2\pi)^{-1}\int _0^{2\pi}M(z+w'e^{i\theta},|u+w_{m+1}e^{i\theta}|)d\theta.
\end{align*}
Thus $M$ is locally plurisubharmonic, and thus subharmonic, in 
$\Omega \times \{ u\in \Complex: |u|>r_0\}$, for some $r_0$.
\end{proof}

Next, we modify $M$ somewhat to obtain a function plurisubharmonic on
$\Omega \times \Complex$.
\begin{lemma}\label{l:M1}
Let $\Omega$, $f$, $M$, and $r_0$ be as in Lemma \ref{l:M}.  For $(z,u)\in
\Omega \times \Complex$, set
$$M_1(z,u)= \left\{ \begin{array}{ll}
M(z,r_0+1)\; & \text{if}\; |u|< r_0+1\\
M(z,|u|)& \text{if}\; |u|\geq r_0+1.
\end{array}\right.
$$
Then $M_1\in \text{PSH}(\Omega \times \Complex)$.
\end{lemma}
\begin{proof}
We again use the fact that being plurisubharmonic is a local property.
Clearly, if $z_0\in \Omega$, $|u_0|\not = r_0+1$, 
then $M_1$ is plurisubharmonic
in a neighborhood of $(z_0,u_0)$.  If $|u_0|=r_0+1$, then, since
$M(z_0,\bullet)$ is increasing and plurisubharmonic, for $z_0\in \Omega$
$$M_1(z_0,u_0)= M(z_0,|u_0|)
\leq (2\pi)^{-1}\int _0^{2\pi}M_1(z_0+w'e^{i\theta},u_0+w_{m+1}e^{i\theta})
d \theta
$$ for all $w=(w',w_{m+1})\in \Complex^{m+1}$, $\|w\|$ sufficiently small.
\end{proof}

With this preparation, we may now prove the following extension of
\cite[Corollary 1.42]{l-g}, 
which we shall apply in Section \ref{s:pot} to prove our
main theorem.
\begin{prop} \label{p:mca}
 Let $\Omega \subset \Complex^m$ be an open, connected set and
let $f$ satisfy assumptions (A0).  Let $\rho(z)$ be the 
order of $g_z(\lambda)=f(z,\lambda)$ in $0<\arg \lambda<\pi$.  If
$\rho(z)\leq \rho_0$ for all $z\in\Omega$, $\rho(z_0)=\rho_0$ for
some $z_0\in \Omega$, and $\rho_0>\max(\alpha,1)$, then $\rho(z)=\rho_0$ for
all $z\in \Omega \setminus E$, where $E\subset \Omega $ is a pluripolar
set.
\end{prop}
\begin{proof}
Choose a $\beta \in \Real$ such that $\max(\alpha, 1)<\beta<\rho_0$, and 
let $M_1(z,u)$ be as defined in Lemma \ref{l:M1}. 
Note that the order of $u\mapsto M_1(z,u)$ is $\max(\rho(z),\beta)$.

Let $\Omega'$ be open, connected, and bounded, with 
$\overline{\Omega'}\subset \Omega$.  Then by 
\cite[Proposition 1.40]{l-g}, there is a sequence $\{\Psi_q\}$ of 
negative plurisubharmonic 
functions on $\Omega'$ such that
$$-(\rho(z))^{-1}=\lim\sup_{q \rightarrow \infty}\Psi_q(z).$$
In addition,
$$\lim \sup _{q\rightarrow \infty}
\left(\Psi_q(z)+\frac{1}{\rho_0}\right)\leq 0,$$
and $$\lim \sup_{q\rightarrow \infty}
\left(\Psi_q(z_0)+\frac{1}{\rho_0}\right)=0.$$  
Thus, by \cite[Proposition 1.39]{l-g}, 
$\rho(z)=\rho_0$ for $z\in \Omega'\setminus E'$, for some pluripolar set 
$E'\subset \Omega'$.
We can cover $\Omega$ with $\Omega'$ having the properties
as above.  The set $E$ is the
union of the corresponding sets $E'$, and is thus pluripolar.
\end{proof}
We note that this proposition could also be proved by adapting arguments of
\cite{kiselman}.
\section{The scattering matrix and its determinant}\label{s:sd}

In this section, we collect some facts about the scattering matrix and 
its determinant.  
For $V\in L^{\infty}_{\comp}(\Real^d;\Complex)$, let $S_V(\lambda)$ be
the associated scattering matrix and let
$$s_V(\lambda)=\det S_V(\lambda).$$ 
It is meromorphic in the upper half-plane,
with at most a finite number of poles there. 
This is a useful function in the 
study of resonances 
because for odd $d$
its zeros in the upper half-plane coincide, with at most a 
finite number of exceptions, with poles of the resolvent in the lower
half-plane, and the multiplicities agree
(see \cite[(3.7)]{zwpf} or \cite{gu-zw}).  That is, for all but finitely 
many $\lambda_0$, if $\Im \lambda_0>0$
is a zero of order $m_0$ of $s_V(\lambda)$, then $-\lambda_0$ is a pole
of order $m_0$ of $\chi R_V(\lambda)\chi$.

Recall that for 
$V\in L^{\infty}_{\comp}(\Real^d;\Complex)$ the scattering matrix
associated to $\Delta +V$ is given by 
\begin{equation}\label{eq:sm}
S_V(\lambda)=I+c_d\lambda^{d-2}\pi_{\lambda}(V-VR_{V}(\lambda)V)
\pi_{-\lambda}^t
\end{equation}
where 
$\pi_{\lambda}$ is given by
$$(\pi_{\lambda}f)(\omega)=\int e^{-i\lambda x\cdot \omega}f(x)dx.$$  Here
$c_d$ is independent of $\lambda$.

\begin{lemma}\label{l:derivest}
Let $V\in L^{\infty}_{\comp}(\Real^d;\Complex)$.
 For $\lambda \in \Real$, there is a $C_{V}$
so that
$$\left| \frac{d}{d\lambda}\log s_{V}(\lambda)\right|
\leq C_{V}|\lambda|^{d-2}$$
whenever $|\lambda|$ is sufficiently large.
\end{lemma}
\begin{proof}
We use 
$$\frac{d}{d\lambda}\log s_{V}(\lambda)= \tr \left( (S_V(\lambda))^{-1}
\frac{d}{d\lambda}S_{V}(\lambda)\right).$$
If $V$ is real-valued, then for $\lambda \in \Real$, $S_V(\lambda)$ is 
unitary.  Otherwise, we will use (\ref{eq:sm}) to bound
$$\|(S_V(\lambda))^{-1}\| = \|S_V(-\lambda)\| $$
when $\lambda \in \Real.$
By (\ref{eq:sm}), 
$$\|S_V(\lambda)-I\|=|\lambda|^{d-2}\|\pi_{\lambda}\chi(V-VR_{V}(\lambda)V)
\chi \pi_{-\lambda}^t\|.$$  
where $\chi\in C_c^{\infty}(\Real^d)$ is one on the support of $V$.
Using \cite[Corollary 3.7]{yafaev}, 
$$\| \pi_{\lambda} \chi\|_{L^2(\Real^d)\rightarrow L^2(\Sphere^{d-1})}
\leq C_{\chi}|\lambda|^{-(d-1)/2},\; 
\| \chi \pi_{-\lambda}^t \|_{L^2(\Sphere^{d-1})\rightarrow L^2(\Real^{d})}
\leq C_{\chi}|\lambda|^{-(d-1)/2}
.$$
For
 $\lambda \in \Real$ and  $|\lambda|$ is sufficiently large
(depending on $\|V\|_{\infty}$ and $\supp V$) $\|V-VR_V(\lambda)V\|\leq C_V$.
Thus we have, for such $\lambda$,
$$\|S_V(\lambda)-I\|\leq C_V |\lambda|^{-1}.$$

Next, we bound 
$$\left \| \frac{d}{d\lambda}S_{V}(\lambda)\right\|_1 =
 \left\|
\frac{d}{d\lambda}\left(I+c_d\lambda^{d-2}\pi_{\lambda}\chi(V-VR_{V}(\lambda)V)\chi
)\pi_{-\lambda}^t\right)\right\|_1.$$
We bound this just as in \cite[Lemma 3.3]{froeseodd}, 
using the fact that
$$\|AB\|_1 \leq \|A\|_2\|B\|_2.$$
When $|\lambda|$ is 
sufficiently large, 
$$\|V-VR_{V}(\lambda)V\| \leq C_{V},\; \|\frac{d}{d\lambda}VR_{V}(\lambda)V
\| \leq C_V.$$
Moreover, the Hilbert-Schmidt norms of $\pi_{\lambda}\chi$, 
$\chi\pi_{-\lambda}^t$ can be estimated using their explicit Schwartz kernels
to see that $$\|\pi_{\lambda}\chi\|_{2}\leq C_{\chi}, \; 
\|\chi\pi_{-\lambda}^t\|_2\leq C_{\chi}.$$
The Hilbert-Schmidt norms of the derivatives of these operators are also 
bounded above by constants, so that
\begin{equation}
\left \| \frac{d}{d\lambda}S_{V}(\lambda)\right\|_1 \leq C_V 
(1+|\lambda|^{d-2}).
\end{equation}
This finishes the proof.
\end{proof}

We shall consider holomorphic families of potentials that satisfy the
following conditions.  These potentials form a somewhat more general
class than those of Theorem \ref{thm:intermediate}.
\begin{assumpt}[{\bf A1}]
Let $\Omega\subset \Complex^m$ be an open set, 
and let $$V=V(z,x)\in \hol(\Omega_z; L^{\infty}_{\comp}(\Real^d_x;\Complex)).$$
That is, $V$ is holomorphic in the $z$ variables and takes its
values in compactly supported potentials.
Moreover, we require that there be a fixed set $K_{\Omega}\subset \Real^d$ such
that $\supp (V(z,\bullet))\subset K_{\Omega}$ for all $z\in \Omega$.  
\end{assumpt}

\begin{prop}\label{p:sdetprop}
Let $\Omega \subset \Complex ^m$ be open and suppose that $V(z,x)
$ satisfies assumptions (A1).  Let $K\subset \Omega$ be a compact set.
Then $s_{V(z)}(\lambda)$
has the following properties:
\begin{description}
\item [a] There is a constant $C_{K,0}\geq 0 $ such that 
$s_{V(z)}(\lambda)$ is holomorphic on $\Omega'\times \{\lambda :\Im \lambda 
>C_{K,0}\}$ for any open $\Omega'\subset K$.
\item  [b] The constant $C_{K,0}$ can be chosen so that for $z_0\in K$, if $\lambda_0$ is
a zero of $s_{V(z_0)}(\lambda)$ with $\Im \lambda_0>C_{K,0}$, then $-\lambda_0$
is a pole of $R_{V(z_0)}(\lambda)$, and the multiplicities coincide.
\item [c]
 For $\Im \lambda >C_{K,0}$, $z\in K$, there is a constant $C$ (depending
on $V$ and $K$) so that  
$$|s_{V(z)}(\lambda)|\leq C e^{C|\lambda|^d}.$$
\item [d]
If $z\in K$, $\Im \lambda =C_1>C_{K,0}$, then there is a constant $C$ (
depending on $V$, $K$, and $C_1$) so that 
$$|s_{V(z)}(\lambda)|\leq C e^{C|\lambda|^{d-2}}.$$
\end{description}
\end{prop}
\begin{proof}

When $z\in K$, $\| V(z,\bullet)\|_{L^{\infty}}$ is bounded, so that there
is a $C_{K,0}$ such that $I+V(z)R_0(\lambda)$ is invertible when $\Im \lambda
> C_{K,0}$, $z\in K$.  Thus 
$R_{V(z)}(\lambda)=R_0(\lambda)(I+VR_0(\lambda))^{-1}$ is holomorphic in
$\Omega' \times \{\lambda :\Im \lambda 
>C_{K,0}\}$, and, using the explicit expression for $S_V$, so are
$S_{V(z)}(\lambda)$ and $s_{V(z)}(\lambda)$.

Using the relation $(S_V(\lambda))^{-1}=S_V(-\lambda)$, we see that
zeros of $S_V(\lambda)$ with $\Im \lambda>C_{K,0}$ correspond to poles
of 
$S_V(\lambda)$ with $\Im \lambda<-C_{K,0}$.  If $S_{V(z)}(\lambda)$ is 
holomorphic in 
$\Omega' \times \{\lambda :\Im \lambda 
>C_{K,0}\}$, then for $z_0\in \Omega'$, the zeros
of $s_{V(z_0)}(\lambda)$ with $\Im \lambda >C_{K,0}$ correspond,               
with multiplicity, to the poles of $R_{V(z_{0})}(-\lambda)$ (e.g. \cite{gu-zw,
zwpf}).

Property (c) follows as in \cite{zwpf} or \cite{froeseodd}, using the 
fact that $\supp V(z,\bullet)$ and $\|V(z,\bullet)\|_{\infty}$ are bounded
when $z\in K$.

To prove the final property, we use the fact that
$$|\det (I+A)|\leq e^{\|A\|_1}.$$
Again as in \cite{froeseodd} and Lemma \ref{l:derivest}, we have that 
$$\|S_V(\lambda)-I\|_{1} 
\leq  C|\lambda|^{d-2}\|\pi_{\lambda} \chi\|_2 \|V-VR_V(\lambda)V\|
\| \chi \pi_{-\lambda}^t\|_2$$
where $\chi \in C_c^{\infty}(\Real^n)$ is one on the support of $V$.
Using the explicit kernels of $\pi_{\lambda}$ we obtain
$$\|\pi_{\lambda}\chi\|_2 \leq e^{C(|\Im \lambda| +1)}.$$
  We obtain a similar estimate
for $\|\chi \pi_{-\lambda}^t\|_2$.
\end{proof}

\section{Proof of Theorem \ref{t:density}}\label{s:pot}

For integers $p\geq 1$, let $G(u;p)$ be the canonical factor 
$$G(u;p)=(1-u)e^{u+u^2/2+...+u^p/p}.$$  

\begin{lemma}\label{l:derivbd}  Let $d'\in \Integers$, $d'\geq 2.$
Suppose the convergence exponent of the sequence $\{\lambda_j\}$ is strictly
less than $d'$ and that $\lambda_j\not \in \Real$ for all
$j$.  Then there is 
an $\epsilon >0$ and a constant $C_{\epsilon}$ such that for $\lambda \in \Real$, 
$$
\left| \int_0^{\lambda}
\frac{d}{dt}\left( \log \prod_{j=1}^{\infty}G(-t/\lambda_j;d'-1)-
\log \prod_{j=1}^{\infty}G(t/\lambda_j;d'-1)\right)dt \right|
\leq C_{\epsilon}(1+|\lambda|^{d'-\epsilon}).$$
\end{lemma}
\begin{proof}
Let
$$n(r)=\#\{ \lambda_j:|\lambda_j|<r\}.$$
Since the convergence exponent for the sequence is less than $d'$, there
is an $\epsilon>0$ and a constant $C_{\epsilon}$ 
so that $n(r)\leq C_{\epsilon}(1+r^{d'-\epsilon}).$

For the real part of the integral, we use
\begin{multline}
\int_0^{\lambda}\Re \frac{d}{dt}\left( \log \prod_{j=1}^{\infty}G(-t/\lambda_j;d'-1)-
\log \prod_{j=1}^{\infty}G(t/\lambda_j;d'-1)\right)dt \\
=  \log \left|\prod_{j=1}^{\infty}G(-\lambda/\lambda_j;d'-1)\right|-
\log \left|\prod_{j=1}^{\infty}G(\lambda/\lambda_j;d'-1)\right|.
\end{multline}
Then applying standard estimates for canonical products (e.g. \cite[Theorem
I.6]{levin}), we have 
$$\left| \log \prod_{j=1}^{\infty}|G(-\lambda/\lambda_j;d'-1)|-
\log \prod_{j=1}^{\infty}|G(\lambda/\lambda_j;d'-1)|\right| \leq C_{\epsilon}
(1+|\lambda|^{d'-\epsilon})$$
for some $\epsilon >0$.

For the imaginary part of the integral, we use
$$|\arg(1-\lambda/\lambda_j)-\arg(1)|\leq \pi$$
where the argument is chosen to be a continuous function of $\lambda$.
We also use
$$\left|\Im \sum_{p=1}^{d'-1}\frac{1}{p}  u^p\right| \leq C 
|u|^{d'-1}\; \text{if}\; |u|\geq 1 $$
and 
$$ | \arg G(u,d'-1)-\arg (1)| \leq C|u|^{d'} \; \text{if}\; |u|\leq 1.$$

Then, using arguments similar to \cite[Lemma 1.3 and Theorem 1.4]{levin},
\begin{align*}& 
\left| \Im \left( \int_0^{\lambda}
\frac{d}{dt}\left( \log \prod_{j=1}^{\infty}G(-t/\lambda_j;d'-1)-
\log \prod_{j=1}^{\infty}G(t/\lambda_j;d'-1)\right)dt \right)\right|\\
& = \left| \arg \prod_{j=1}^{\infty}G(-t/\lambda_j;d'-1)-
\arg \prod_{j=1}^{\infty}G(t/\lambda_j;d'-1)\right|\\
& \leq Cn(2\lambda)+ |\lambda|^{d'-1} \int _1^{2\lambda} t^{1-d'}dn(t)
+
|\lambda|^{d'} \int _{2\lambda}^{\infty}
t^{-d'}dn(t) \\ &
\leq C_{\epsilon}(1+|\lambda|^{d'-\epsilon}).
\end{align*}
\end{proof}

\begin{lemma}\label{l:cexp}
Let $d\geq 3$ be odd, and $V\in L^{\infty}_{\comp}(\Real^d;\Complex)$.  Then
$s_V(\lambda)$ is of order $d$ in the half-plane $\{\lambda\in \Complex:
\Im \lambda > 2\|V\|_{\infty}+1\}$ if and only if $\mcr_V$ has convergence
exponent $d$.
\end{lemma}
\begin{proof}
We remark that since $N_V(r)\leq C(r^d+1)$, the convergence exponent of
$\mcr_V$ is at most $d$.

We shall actually prove the contrapositive of this lemma.

Let $g(\lambda)$ be holomorphic in a neighborhood of the closed upper 
half-plane, and let $n_g(r)$ be the number of zeros of $g$ in the upper 
half-plane with norm less than $r$, counted with multiplicity.  By 
using intermediate steps from the proof of \cite[Lemma 3.2]{froeseodd}, 
$$\int_0^r\frac{n_g(t)}{t}dt =\frac{1}{2\pi}
\int_0^r t^{-1}\int_{-t}^t\frac{g'(s)}{g(s)}ds dt +\frac{1}{2\pi}\int_0^{\pi}
\log |g(re^{i\theta})|d\theta.$$
We apply this to $s_V(\lambda)$ (multiplied by a suitable polynomial if
it has poles in the upper half-plane).  Using Lemma \ref{l:derivest} as
well, we 
see that if
 $s_V$ has order strictly less than $d$ in this region,
$\mcr_V$ has convergence exponent strictly less than
$d$.

Now suppose that $\mcr_V$ has convergence exponent
$\rho$ strictly less than $d$.
We may write \cite{zwpf}
$$s_V(\lambda)= \alpha e^{ig(\lambda)}
\frac{P(-\lambda)}{P(\lambda)}$$
where $\alpha$ is a constant,
$$P(\lambda)=\prod_{\lambda_j\in \mcr_V, \lambda_j \not = 0}
G(\lambda/\lambda_j;d-1)
$$
and $g(\lambda)$ is a polynomial of order at most $d$.
The canonical product $P(\lambda)$ is of order $\max(\rho,d-1)$.  By
the minimum modulus theorem, then,
$P(\lambda)/P(-\lambda)$ is of order $
\max(\rho,d-1)$ in the upper half plane in 
question.  From Lemma 
\ref{l:derivbd}, we know that 
$$
\left| \int_0^{\lambda}\frac{d}{dt}(\log P(t)-\log P(-t))dt\right|
\leq C_{\epsilon}(1+|\lambda|^{d-\epsilon})$$
for some $\epsilon>0$.  Thus, using Lemma
\ref{l:derivest}, we see that $g$ must
have order less than $d$.  Therefore, $s_V(\lambda)$ has order strictly
less than $d$ in 
$\{\lambda\in \Complex:
\Im \lambda \geq 2\|V\|_{\infty}+1\}.$
\end{proof}

Our main theorem allows a more general family of potentials
than Theorem \ref{thm:intermediate}.
\begin{thm}\label{t:major1}
Let $d\geq 3$ be odd and let $\Omega \subset \Complex^m$ be open 
and connected.
Suppose $V(z,x)$ satisfies assumptions (A1), and for some $z_0\in \Omega$,
$\mcr_{V(z_0)}$ has convergence exponent $d$.  Then $\mcr_{V(z)}$ has
convergence exponent $d$ for all $z\in \Omega \setminus E$, where
$E\subset \Omega$ is a pluripolar set.
\end{thm}
\begin{proof}  By Lemma \ref{l:cexp}, the order of $s_{V(z_0)}(\lambda)$ 
in the 
upper half-plane is $d$.  Given an open, connected $\Omega'\subset \Omega$
such that $\overline{\Omega}'\subset \Omega $ is bounded, 
by Proposition \ref{p:sdetprop} we may
apply Proposition \ref{p:mca} to 
$s_{V(z)}(\lambda+iC_{\overline{\Omega'}}+i)$.
From Proposition \ref{p:mca} we see that in the upper half-plane
$s_{V(z)}(\lambda +
iC_{\overline{\Omega'}}+i
)$ 
 has order $d$ for $z\in \Omega' \setminus E'$, for some pluripolar set 
$E'$.  Again by Lemma \ref{l:cexp}, this means that the convergence
exponent of $\mcr_{V(z)}$ is $d$ for $z\in \Omega' \setminus E'$.  As in
the proof of Proposition \ref{p:mca}, we can cover $\Omega$ by such 
$\Omega'$, and the set $E$ is the union of the corresponding sets $E'$. 
\end{proof}

\section {A class of potentials with fixed sign}\label{s:fixedsign}

For $V\in L^{\infty}_{\comp}(\Real^d;\Complex)$ let 
\begin{equation}
B_V(\lambda)=\frac{V}{|V|^{1/2}}R_0(\lambda)|V|^{1/2}.
\end{equation} 
Then the poles of $R_V(\lambda)$ are the 
zeros of $I+B_V(\lambda)$.  In this section we use this fact and a study 
of related holomorphic functions to prove Theorem \ref{thm:fixedsign}.
Throughout this section we assume that $d$ is odd.

\subsection{Lower bounds on a determinant}\label{s:lbd}
In this subsection we obtain lower bounds on 
$\det(I+B_V^{2m}(\lambda))$ when $m>d/4$ and $V$ can be bounded below by the 
characteristic function of the ball.  In the next subsection we will use 
this lower bound and some results of several complex variables to 
prove Theorem \ref{thm:fixedsign}.

\begin{lemma}\label{l:detcomparison}
 Let $V\in L^{\infty}_{\comp}(\Real^d;\Real)$ satisfy
$V\geq \chi_{B(a,x_0)}$ where $\chi_{B(a,x_0)}$ is the characteristic
function of the ball of radius $a>0$ centered at $x_0$.  Let 
$V_0=\chi_{B(a,0)}$.  Then, for $s\in \Real_+$, $m>d/4$, $m\in \Integers$,
$$\det(I+(B_V(-is))^{2m})\geq 
\det(I+(B_{V_0}(-is))^{2m}).$$
\end{lemma}
Before proving the lemma, we remark that the sign in front of 
$(B_{V}(-is))^{2m}$ may appear puzzling,
as the zeros of $\det(I+(B_V(\lambda))^{2m})$ do not, in general,
include the poles of $R_V$, 
while those of $\det(I-(B_V(\lambda))^{2m})$
do (compare Lemma \ref{l:zeros}).  The 
sign is positive
so that we may work with the determinant of the identity plus a positive 
operator.
In the proof of Proposition \ref{p:sumorder} we introduce a complex 
parameter, and this allows us to reconcile the apparent differences.
\begin{proof}
If $A$ is a trace class operator, 
\begin{equation}\label{eq:determinant}
\det(I+A)=\prod_j(I+\mu_j(A)),
\end{equation}
 where $\mu_j$ are the eigenvalues of
$A$ repeated according to their multiplicity, and 
$|\mu_1(A)|\geq |\mu_2(A)|\geq.... $ 

We note that for $s\in \Real_+$, $B_V^2(-is)$ is a 
positive, self-adjoint operator, so that all of its eigenvalues are
non-negative.  Let $V_1=\chi_{B(a,x_0)}$.  Then, using the max-min principle,
$$\mu_j(B_V^2(-is))\geq \mu_j(B_{V_1}^2(-is))\geq 0$$
and thus $$
\mu_j(B_V^{2m}(-is))\geq \mu_j(B_{V_1}^{2m}(-is))\geq 0.$$  Since
the eigenvalues of $B_{V_1}^{2m}$ are the same as those of $B_{V_0}^{2m}$,
using (\ref{eq:determinant}) finishes the proof of the lemma.
\end{proof}

Next we will describe the resolvent $R_0(\lambda)$ in a way which will be
useful for our purposes.  
Let $\sigma_1^2< \sigma_2^2\leq ...$ be the eigenvalues
of the Laplacian on the sphere $\Sphere^{d-1}$, repeated according to their
multiplicity, and let $\{\phi_j\}$ be a corresponding set of orthonormal
eigenfunctions.  We use the notation of \cite{olver} for the Bessel functions
$J_{\nu}$, the modified Bessel functions $I_{\nu}$ and
$K_{\nu}$, and the Hankel function $H^{(1)}_{\nu}$.  Let $\lambda \in 
\Complex$ and $(r,y)\in \Real_+\times \Sphere^{d-1}$ be polar coordinates
on $\Real^d$.  Then 
\begin{align}\label{eq:R0}
&\frac{2i}{\pi}(R_0(\lambda)f)(r,y)\\ 
\nonumber 
&= \sum_{k=1}^{\infty}\sum_{\sigma_j^2=k(k+d-2)}
 \int_0^r\int_{\Sphere^{d-1}}\frac
{H^{(1)}_{\nu_k}(\lambda r)J_{\nu_k}(\lambda r')}{(rr')^{(n-2)/2}} 
\phi_j(y)\overline{\phi}_j(y')f(r',y')(r')^{n-1}d\sigma_{y'}dr'\\ 
\nonumber
&\mbox{  }
+ \sum_{k=1}^{\infty}\sum_{\sigma_j^2=k(k+d-2)}
 \int_r^{\infty}\int_{\Sphere^{d-1}}
\frac{H^{(1)}_{\nu_k}(\lambda r')J_{\nu_k}(\lambda r)}{(rr')^{(n-2)/2}}
 \phi_j(y)\overline{\phi}_j(y')f(r',y')(r')^{n-1}d\sigma_{y'}dr'.
\end{align}
Here 
$$\nu_k=k+\frac{d}{2}-1.$$

To obtain a lower bound on the eigenvalues of $B_{V_0}^2(-is)$ we shall 
need some lower bounds on Bessel functions.
\begin{lemma}\label{l:besselbound}
Let $s,\;M\in \Real_+$ and $\nu -1/2 \in \Natural.$  
Then there is a constant $c>0$ such that
\begin{align*}
|J_{\nu}(-i\nu s)|& \geq c\frac{e^{c\nu}}{\sqrt{\nu}}\\
|H_{\nu}^{(1)}(-i\nu s)|& \geq c\frac{e^{c\nu}}{\sqrt{\nu}}
\end{align*}
when $3<s<M$ and $\nu$ is sufficiently large.
\end{lemma}
\begin{proof} We use, from \cite[9.6.3 and 9.6.30]{olverbig}
\begin{equation}\label{eq:j}
|J_{\nu}(-is)|=|I_{\nu}(s)|
\end{equation}
and, from \cite[9.1.40, 9.6.4, and 9.6.31]{olverbig}
\begin{equation}\label{eq:H}
H^{(1)}_{\nu}(-is)=\frac{-2i}{\pi}e^{-3\nu \pi i/2}K_{\nu}(s)
- 2e^{-\nu \pi i/2} I_{\nu}(s).
\end{equation}
For $3\leq z \leq M<\infty$ there are  constants $c,\; C>0$ such that
\begin{align*}
I_{\nu}(\nu z)& \geq c \frac{e^{\nu \xi}}{\sqrt{\nu}}\\
|K_{\nu}(\nu z)|&\leq C e^{-\nu\xi}
\end{align*}
when $\nu$ is sufficiently large and 
$\xi =(1+z^2)^{1/2}+\ln \frac{z}{1+(1+z^2)^{1/2}}$ \cite[10.7.16]{olver}.  
Applying this and (\ref{eq:j}), for some
$c>0$
$$|J_{\nu}(-i\nu s)|\geq c \frac{e^{c\nu}}{\sqrt{\nu}}$$
when $3\leq s\leq M$ and $\nu$ is sufficiently large.  
Similarly, using
(\ref{eq:H}), the upper bound on $K_{\nu}(\nu z)$ and the lower bound on
$I_{\nu}(\nu z)$, we obtain the second part of the lemma.
\end{proof}

The following lemma shows
that the holomorphic 
function  $\det(I+(B_{V_0}(-is))^{2m})$ 
introduced in Lemma \ref{l:detcomparison} has order at least $d$.
\begin{lemma}\label{l:orderd}
For $a>0$, let 
$V_0=\chi_{B(a,0)}\in L^{\infty}_{\comp}(\Real^d)$.  
Then, for $s\in \Real_+$, $m>d/4$, $m\in \Integers$, there is a constant
$c>0$ such that
$$
\det(I+(B_{V_0}(-is))^{2m})\geq c e^{cs^d}$$
when $s$ is sufficiently large.
\end{lemma}
\begin{proof}
We first obtain a lower bound on some of the eigenvalues of 
$B^2_{V_0}(-is)=(V_0^{1/2}R_0(-is)|V_0|^{1/2})^2$.  
Since $V_0$ is radial, we can 
write  $B^2_{V_0}(-is)=\sum_{k}B^2_{k,V_0}(-is)$ where 
$B_{k,V_0}(-is)$ acts on the eigenspace of the Laplacian on $\Sphere^{d-1}$ 
with eigenvalue $k(k+d-2)$ and multiplicity  $m(k)\geq ck^{d-2}$,
for some $c>0$.  We will bound $\|B_{k,V_0}(-is)\|$ from below, giving 
us a lower bound on $m(k)$ of the eigenvalues of $B^2_{k,V_0}(-is)$.

Using (\ref{eq:R0}), 
\begin{align*}
\|B_{V_0}(-is)\|\geq 
\|\chi_{[0,a/2]}r^{-(n-2)/2}J_{\nu_k}(-isr)\|_{L^2(\Real^d)}
\|\chi_{[a/2,a]}r^{-(n-2)/2}H^{(1)}_{\nu_k}(-isr)\|_{L^2(\Real^d)}
\end{align*}
where $\chi_{[\alpha, \beta]}$ is the characteristic function of the 
interval $[\alpha,\beta]$.  
By Lemma \ref{l:besselbound}, for 
$M\in \Real_+$ and $\frac{sa}{2M}
<\nu_k<\frac{sa}{12}$ there is some $c>0$ such that
\begin{align}
\|\chi_{[0,a/2]}r^{-(n-2)/2}J_{\nu_k}(-isr)\|_{L^2(\Real^d)}^2 &
\geq \int_{a/4}^{a/2} c r \frac{e^{c\nu_k}}{\nu_k}dr\\ 
\nonumber
& \geq \frac{c}{\nu_k}e^{c\nu_k}.
\end{align}
Here and throughout $c$ is a positive constant whose value may change
from line to line.  Using Lemma \ref{l:besselbound} in a similar way, for
$M'\in \Real_+$,
\begin{equation}
\|\chi_{[a/2,a]}r^{-(n-2)/2}H_{\nu_k}^{(1)}(-isr)\|^2_{L^2(\Real^d)}
\geq c \frac{e^{c\nu_k}}{\nu_k}
\end{equation}
for some $c>0$ when $\frac{sa}{M'}<\nu_k<\frac{sa}{6}$.  
Thus, with $\alpha =d/2-1$, 
$$\|B_{V_0,k}(-is)\|\geq \frac{c}{k+\alpha}e^{ck}$$
when $\frac{sa}{M'}<k+\alpha<\frac{sa}{12}.$  Thus $m_k\geq c k^{d-2}$
eigenvalues of $B^2_{V_0,k}$ are at least as large as $c (k+\alpha)^{-2}e^{ck}$
when $\frac{sa}{M'}<k+\alpha<\frac{sa}{12}.$
Then, taking $M'$ and $s$ sufficiently large,
\begin{align*}
\det(I+(B_{V_0}(-is))^{2m}) & 
\geq  \prod_{\frac{sa}{M'}<k+\alpha<\frac{sa}{12}}
\left(1+ c\frac{e^{ck}}{(k+\alpha)^{2m}}\right)^{ck^{d-2}}\\ & 
\geq \exp \sum_{\frac{sa}{M'}<k+\alpha<\frac{sa}{12}}
(ck^{d-2}-ck^{d-2}\ln (k+\alpha)+ck^{d-1})\\
& \geq c \exp(cs^d).
\end{align*}
\end{proof}

\subsection{Proof of Theorem \ref{thm:fixedsign} and Corollary 
\ref{t:density} }
In this subsection we use the results of Section \ref{s:lbd}, Theorem
\ref{thm:intermediate}, and some results from \cite{l-g} to prove Theorem
\ref{thm:fixedsign}.  We also prove Corollary \ref{t:density}.

If $m>d/2$ is an integer, $\det(I-(-1)^mB_V^m(\lambda))$ is a holomorphic
function of $\lambda$.  Moreover, its zeros include the 
poles of $R_{V}(\lambda)$.  In fact, we can say more, as the next
lemma shows (compare \cite[Proposition 1]{zwodd}).
\begin{lemma}\label{l:zeros}
 Let $m>d/2$ be an integer and let $\omega=e^{2\pi i/m}$.
Let $V\in L^{\infty}_{\comp}(\Real^d;\Complex)$.
Then the zeros of $\det(I-(-1)^mB_V^m(\lambda))$ correspond, with multiplicity,
to $\cup_{k=1}^m \mcr_{\omega^k V}$.
\end{lemma}
\begin{proof}
We note that
$$I-(-1)^mB_V^m(\lambda)=\prod_{k=1}^m
(I+\omega^k B_V(\lambda))=\prod_{k=1}^m(I+B_{\omega^kV}(\lambda)).$$  
The lemma follows from using the
fact that the zeros of $I+B_V$ correspond, with
multiplicity, to
$\mcr_V$.
\end{proof}

We shall need some knowledge of $\det(I-(-1)^mB_V^m(\lambda))$ in the 
upper half plane.  This is analogous to \cite[Lemma 3.2]{froese}.
\begin{lemma} \label{l:near1}
 Let $V\in L^\infty_{\comp}(\Real^d;\Complex),$
 $m>d/2$ be an integer and 
$$h_V(\lambda) =\det(I-(-1)^mB_V^m(\lambda)).$$
Then  for $0<\theta<\pi$, $\epsilon >0$, there is a 
$C$ (depending on $V$, $\theta$, and $\epsilon$ )
such that for $r\in \Real_+$, 
$$|h_V(re^{i\theta})-1|\leq C r^{\epsilon -1}.$$
\end{lemma}
\begin{proof}
We use the fact that $|\det(I+A)-1|\leq \|A\|_1 e^{\|A\|_1+1}$ 
\cite[Theorem XIII.104]{r-siv}.  For 
$\chi\in C_c^{\infty}(\Real^d)$, 
$$\|\chi R_0(re^{i\theta})\chi\|_{H^s(\Real^d)\rightarrow H^{s+2}(\Real^d)}\leq C$$ and
$$\|\chi R_0(re^{i\theta})\chi\|_{H^s(\Real^d)\rightarrow
 H^{s}(\Real^d)}\leq 
\frac{C}{r^2}.$$  Here $C$ denotes a positive constant whose value changes from
line to line and may depend on parameters other than $r$.  Therefore, 
for any $\epsilon'>0$ and $p>d/2$, 
$$\| \chi R_0(re^{i\theta})\chi\|_p\leq \frac{C}{r^{2-d/p-\epsilon'}}$$
and $$\|B_V^m(re^{i\theta})\|_1\leq \frac{C}{r^{1-\epsilon}}.$$
\end{proof}

\begin{prop}\label{p:sumorder}
Let $V\in L^{\infty}_{\comp}(\Real^d;\Real)$ be bounded below by the
characteristic function of a ball.  Let $m>d/4$ be an 
integer and let $\omega =e^{\frac{\pi i}{m}}.$  Then
$$\lim \sup_{r\rightarrow \infty}\frac{\log (\sum_{j=1}^{2m}N_{\omega^jzV}(r))}
{\log r}=d$$
for all $z\in \Complex \setminus E$, where $E$ is a pluripolar set.
\end{prop}
\begin{proof}
Consider the function $\det(I-B^{2m}_{zV}(\lambda))$.  This is a holomorphic
function of $(\lambda, z)\in \Complex^2.$  Moreover, as in 
\cite[Proposition 3]{zwodd}, it is (for fixed
$z$) of order at most $d$ in $\lambda$.  On the other hand, if $z^{2m}=-1$,
by Lemma \ref{l:orderd} it is of order at least $d$ in $\lambda$.  Thus, 
applying \cite[Propositions 1.39 and 1.40]{l-g}
as in the proof of Proposition \ref{p:mca}, $\det(I-B^{2m}_{zV}(\lambda))$
is of order $d$ for $z\in \Complex \setminus E$ for a pluripolar set $E$.

Now fix $z\in \Complex \setminus E$.  Suppose 
$$\lim \sup_{r\rightarrow \infty}
\frac{\log (\sum_{j=1}^{2m}N_{\omega^j z V}(r))}
{\log r}=d'<d.$$  Then we may write
\begin{equation}\label{eq:detfactor}
\det(I-B^{2m}_{zV}(\lambda))=\alpha_z e^{ig_z(\lambda)}
\prod_{\begin{subarray}{c}
\lambda_j\in \cup_k \mcr_{\omega^k z V},\ \\\lambda_j\not =0
\end{subarray}}
G(\lambda/\lambda_j;d-1)
\end{equation}
where $\alpha_z$ is a constant,
$g_z$ is a polynomial of order at most $d$, and $G(\zeta;d-1)$ is
the canonical factor of order $d-1$.  There are at most finitely many
elements of $\cup_k \mcr_{\omega ^k z V}$ in the 
upper half plane.  Thus standard estimates on canonical products 
and the minimum modulus theorem show that
 for $\Im \lambda$ sufficiently large
and $0<\theta_1<\arg \lambda<\theta_2<\pi$,
the canonical product in (\ref{eq:detfactor}) must satisfy, for every
$\epsilon >0$ and some $C_{\epsilon}$,
$$\frac{1}{C_{\epsilon}}e^{-C_{\epsilon}|\lambda|^{d'+\epsilon}}
\leq \left| \prod_{\begin{subarray}{c}\lambda_j\in \cup_k \mcr_{\omega^k z V}
\\ \lambda_j \not = 0 \end{subarray}}
G(\lambda/\lambda_j;d-1)\right|
\leq C_{\epsilon}e^{C_{\epsilon}|\lambda|^{d'+\epsilon}}.$$
Thus, using Lemma \ref{l:near1}, $g_z(\lambda)$ must be of order strictly 
less than $d$, and so $\det(I-B^{2m}_{zV}(\lambda))$ is of order strictly
less than $d$ in $\lambda$, a contradiction.
\end{proof}

Now we can give the proof of Theorem \ref{thm:fixedsign}.  
\begin{proof}[Proof of Theorem \ref{thm:fixedsign}]
Fix an integer $m$ with $m>d/4$.
By Proposition \ref{p:sumorder} and using the notation of that
proposition, 
$$\lim \sup_{r\rightarrow \infty}\frac{\log (\sum_{j=1}^{2m}N_{\omega^jzV}(r))}
{\log r}=d$$
for $z\in \Complex \setminus E'$, some pluripolar set $E'$.  If
$z_1\in \Complex \setminus E'$, then, there is some $j_1$ such that
$$\lim \sup _{r\rightarrow \infty}\frac{\log N_{\omega^{j_1}z_1V}(r)}
{\log r}=d.$$
Thus by Lemma \ref{l:cexp}
the potential $V_1(z,x)=zV(x)$ satisfies the assumptions of Theorem 
\ref{thm:intermediate} with $z_0=\omega^{j_1}z_1$.  Applying 
Theorem \ref{thm:intermediate} finishes the proof. 
\end{proof}

We may now give the proof of Corollary \ref{t:density}.
\begin{proof}[Proof of Corollary \ref{t:density}] 
 Let $V_0\in L^{\infty}_{\comp}
(\Real^d;\Complex)$ and let $\epsilon>0$.  Let 
$V_1\in C_c^{\infty}(\Real^d;\Real)$ be bounded below by the characteristic 
function of a ball.  Then by Theorem
\ref{thm:fixedsign} $\mcr_{zV_1}$ has convergence
exponent $d$ for all $z\in \Complex \setminus E_1$ for some pluripolar set
$E_1$.  The set $E_1\subset \Complex\simeq \Real^2$ not only has Lebesgue 
measure
zero, but 
its restriction to 
$\Real$ is of Lebesgue measure zero in $\Real$
(e.g. \cite[Section 3.2]{ransford}).  Thus we may choose $z_1\in 
\Real \setminus (E_1)\restrict_{ \Real}$ so that 
$\|z_1V_1\|_{L^{\infty}}<\epsilon/2$
and $\mcr_{z_1V_1}$ has convergence exponent $d$.


Now consider $V(z)=z z_1V_1+(1-z)V_0$.  Then $V(z,x)$ is in the 
framework of Theorem \ref{thm:intermediate}, 
$V(1)=z_1V_1$, and $V(0)=V_0$.  
Using Theorem
\ref{thm:intermediate}, $\mcr_{V(z)}$ has convergence exponent $d$ for all 
$z\in \Complex \setminus E$, where $E$ is a pluripolar set.  
 Thus we may find a 
$z_2\in \Real\setminus (E\restrict _{\Real})$ with $|z_2|<\epsilon (2
(\|V_0\|_{L^{\infty}}+1))^{-1}$.
The potential $V_2(x)=V(z_2,x)$ thus has $\mcr_{V_2}$ with convergence exponent
$d$ and $\| V_0-V_2\|_{L^{\infty}}<\epsilon$.  We remark that if $V_0$ is
real-valued, then so is $V_2$.

With fairly straightforward modifications, the same proof gives the result for
smooth potentials in the $C^{\infty}$ topology.
\end{proof}

\small
\noindent
{\sc 
Department of Mathematics\\
University of Missouri\\
Columbia, Missouri 65211\\}
\end{document}